\documentclass[11pt,leqno]{amsart}
\usepackage[T1]{fontenc}
\usepackage[a4paper,left=25mm,right=25mm,top=25mm,bottom=28mm,
  headheight=14pt,headsep=7mm,footskip=12mm]{geometry}
\usepackage{amssymb,mathtools}
\usepackage{microtype}
\usepackage[hidelinks,bookmarksnumbered]{hyperref}
\hypersetup{
  pdftitle={Gaussian limits for non-cyclic and elementary abelian Sylow subgroups of the multiplicative group},
  pdfauthor={Aimin Guo},
  pdfkeywords={Sylow subgroups, elementary abelian groups,
    joint Gaussian limits, Carmichael's function,
    common prime divisors, arithmetic covariance}
}
\theoremstyle{plain}
\newtheorem{theorem}{Theorem}[section]
\newtheorem{lemma}[theorem]{Lemma}
\newtheorem{corollary}[theorem]{Corollary}
\theoremstyle{remark}
\newtheorem{remark}[theorem]{Remark}
\allowdisplaybreaks[4]
\newcommand{\E}{\mathbb E}
\newcommand{\Pp}{\mathbb P}
\newcommand{\ind}{\mathbf 1}
\newcommand{\Var}{\operatorname{Var}}
\newcommand{\Cov}{\operatorname{Cov}}
\newcommand{\rad}{\operatorname{rad}}
\newcommand{\dd}{\,\mathrm d}
\title[Gaussian limits for Sylow subgroup counts]
{Gaussian limits for non-cyclic and elementary abelian\\
Sylow subgroups of the multiplicative group}
\author[A. Guo]{Aimin Guo}
\address{School of Mathematics and Statistics, Anhui Normal University,
Wuhu 241002, Anhui, P. R. China}
\email{18339681864@163.com}
\date{}
\begin{document}
\begin{abstract}
We study the number of non-cyclic Sylow subgroups
of $(\mathbb Z/n\mathbb Z)^\times$. Pollack suggested investigating
whether this count has a limiting distribution after suitable
normalization. For uniformly chosen $n\le x$, we prove a Gaussian
limit with variance asymptotic to
$\tfrac12\log_2x/\log_3x$. Motivated by the fixed-rank normal-order
results of Morales, Polanco, and Pollack, we further prove joint
Gaussian limits for this count and the numbers of elementary
abelian Sylow subgroups of any fixed finite collection of ranks
at least two. We give explicit prime-sum centers and
asymptotics for all arithmetic covariances. The same joint theorem
includes the number of common prime divisors of Euler's and
Dedekind's functions. The proof uses the integer Kubilius
model, estimates for the dependence between distinct output
primes, and mean-square bounds for arithmetic truncation.
\end{abstract}
\maketitle

\section{Introduction}
Let
\[
 W(n)=\omega\left(\frac{\varphi(n)}{\lambda(n)}\right).
\]
Here $\omega$ counts distinct prime divisors, $\varphi$ is Euler's
function, and $\lambda$ is Carmichael's function, the exponent of
$(\mathbb Z/n\mathbb Z)^\times$. Since this group is abelian,
$W(n)$ is the number of its non-cyclic Sylow subgroups.

In 2021, Pollack \cite[Theorem 1]{Pollack} proved that $W(n)$
has normal order $\log_2n/\log_3n$. This gives the typical
size of $W(n)$. It does not determine the size of its fluctuations,
or the proportion of integers for which $W(n)$ is larger or smaller
than its center. In his concluding remarks
\cite[\S4]{Pollack}, he recalled Gaussian laws for
$\omega(\varphi(n))$ and logarithmic subgroup counts. He then
suggested investigating whether $W(n)$ also has a limiting
distribution after suitable normalization.

Our first result gives such a limit. We prove that
$W(n)$ has a Gaussian limit with
an explicit center, and its arithmetic variance is asymptotic
to $\tfrac12\log_2x/\log_3x$; see Corollary \ref{cor:scalar}.
This determines the limiting proportions of integers $n\le x$
for which the standardized count lies in any fixed interval.

In 2026, Morales, Polanco, and Pollack
\cite[Theorem 1.1]{MPP} proved a normal-order theorem for
elementary abelian Sylow subgroups of fixed rank. For $r\ge2$,
let $E_r(n)$ count the primes $q$ for which the entire Sylow
$q$-subgroup of $(\mathbb Z/n\mathbb Z)^\times$ is isomorphic to
$(\mathbb Z/q\mathbb Z)^r$.
They showed that $E_r(n)$ has normal order
$\log_2n/(r(r-1)\log_3n)$. This gives the typical size of each count.
We now ask how these counts fluctuate and how they vary
together with $W(n)$.
Theorem \ref{thm:joint} gives a joint Gaussian
limit for every fixed finite collection of ranks. It also gives
explicit asymptotics for all arithmetic variances and covariances.

As a further application, the joint theorem includes
\[
 C(n)=\omega(\gcd(\varphi(n),\psi(n))),\qquad
 \psi(n)=n\prod_{p\mid n}(1+1/p).
\]
Here $\psi$ is Dedekind's function. Thus $C(n)$ counts the
distinct primes dividing both
$\varphi(n)$ and $\psi(n)$, each prime counted once.
Corollary \ref{cor:sigma} gives the corresponding
sum-of-divisors variant.

Related work studies multiplicities of fixed invariant-factor
orders \cite{MS}, subgroup counts \cite{MTgroups}, and polynomials
in additive functions \cite{MTpoly}.

To study fluctuations, we need smaller error terms than those
used for normal orders. The bound for the mean truncation error in
\cite[Lemma 2.4]{MPP} is $O(\log_2x/(\log_3x)^2)$.
This estimate alone does not show that the truncation is
negligible on the fluctuation scale $\sqrt{\log_2x/\log_3x}$.
We prove that the truncation error has both
$L^1$ and $L^2$ norms smaller than this scale by a factor tending
to zero. The first bound transfers the Gaussian limit; the second
also transfers the variances and covariances in Theorem
\ref{thm:joint}. We also bound the dependence caused by input
primes that contribute to several output primes
$q\asymp\log_2x$.
The integer Kubilius model and an association argument
\cite{Newman} then yield the joint limit. These covariance
asymptotics describe fluctuations that normal-order statements
alone do not determine. The fixed-rank covariance entries
in \eqref{eq:covmatrix} are classical in occupancy theory
\cite[\S5]{BG}; the arithmetic estimates justify their occurrence
here.

All logarithms are natural. We write $\log_1t=\log t$ and
$\log_{j+1}t=\log(\log_jt)$ whenever these expressions are defined.
We write $\Pp_x$ for the uniform probability measure
on $1\le n\le\lfloor x\rfloor$, with expectation $\E_x$,
variance $\Var_x$, and covariance $\Cov_x$.

We use $\ind_A$ for the indicator of an event $A$, and
$\Longrightarrow$ for convergence in distribution.
The notation $N_d(0,\Sigma)$ denotes the centered $d$-dimensional
Gaussian law with covariance matrix $\Sigma$.
For arithmetic random variables, $\|Y\|_2=(\E_x|Y|^2)^{1/2}$.
We use $O$ and $\ll$ for bounds up to a constant, and $o$ for
a ratio tending to zero. Subscripts indicate permitted dependence
on fixed parameters. The positive constants $C_0$ and $C_K$
may change from line to line; $C_0$ is absolute and $C_K$
may depend on $K$.

The letters $p,q,\ell$ denote primes; $r,s,k,K$ denote integer ranks
or indices unless explicitly stated otherwise.
The centers in our limit theorems are the explicit prime sums
\begin{align}
 A_W(x)&=\sum_{q\ge5}
   \left[1-\left(1+\frac{\log_2x}{q-1}\right)e^{-\log_2x/(q-1)}\right],\notag\\
 A_C(x)&=\sum_{q\ge5}\left(1-e^{-\log_2x/(q-1)}\right)^2,\notag\\
 B_r(x)&=\frac1{r!}\sum_{q\ge5}
           \left(\frac{\log_2x}{q-1}\right)^r e^{-\log_2x/(q-1)}
           \qquad(r\ge2).
 \label{eq:centers}
\end{align}
For fixed $x$, the summands are $O_{x,r}(q^{-2})$, so the sums converge.
Every sum indexed by $q$ below is over primes.
The restriction $q\ge5$ removes only two possible output primes.
The primes $2,3$ contribute at most two to each of $W,C,E_r$.
Their bounded contribution is included in the comparison errors
below; it does not affect the $\sqrt{\frac{\log_2x}{\log_3x}}$-scale limit or the leading
$\frac{\log_2x}{\log_3x}$-scale covariances. See Lemma~\ref{lem:groups} and the final
covariance comparison.

\begin{theorem}\label{thm:joint}
Fix an integer $K\ge2$. As $x\to\infty$,
\begin{align}
 &\sqrt{\frac{\log_3x}{\log_2x}}\,
 \Bigl(W(n)-A_W(x),\ C(n)-A_C(x),E_2(n)-B_2(x),\ldots,E_K(n)-B_K(x)\Bigr)\notag\\
 &\hspace{18mm}\Longrightarrow N_{K+1}(0,\Sigma_K)
       \quad\hbox{under }\Pp_x.
 \label{eq:joint}
\end{align}
The indices $0,1,2,\ldots,K$ correspond to
$W,C,E_2,\ldots,E_K$, respectively.
The positive definite covariance matrix is specified by symmetry and
\begin{align}
 (\Sigma_K)_{00}&=\frac12,\qquad
 (\Sigma_K)_{11}=6\log(9/8),\qquad
 (\Sigma_K)_{01}=\log(9/8),\notag\\
 (\Sigma_K)_{0r}&=\frac{r+1}{2^r r(r-1)},\qquad
 (\Sigma_K)_{1r}=\frac{2^{1-r}-3^{1-r}}{r(r-1)}
       \quad(2\le r\le K),\notag\\
 (\Sigma_K)_{rs}&=\frac{\ind_{\{r=s\}}}{r(r-1)}
       -\frac{(r+s-2)!}{2^{r+s-1}r!s!}
       \quad(2\le r,s\le K).
 \label{eq:covmatrix}
\end{align}
The arithmetic means satisfy
$\E_xW=A_W(x)+O_K(\log_3x)$,
$\E_xC=A_C(x)+O_K(\log_3x)$, and
$\E_xE_r=B_r(x)+O_K(\log_3x)$ for $2\le r\le K$.
Since $\log_3x=o(\sqrt{\frac{\log_2x}{\log_3x}})$, these mean estimates imply that
the centers in \eqref{eq:joint} can simultaneously be replaced
by the exact arithmetic means.

Moreover, the arithmetic covariance matrix has the entrywise asymptotic
\[
 \frac{\log_3x}{\log_2x}
       \Cov_x\bigl((W,C,E_2,\ldots,E_K)\bigr)
       =\Sigma_K+o_K(1).
\]
The centers have the leading asymptotics
\begin{align*}
 A_W(x)&\sim\frac{\log_2x}{\log_3x},\qquad
 A_C(x)\sim(2\log2)\frac{\log_2x}{\log_3x},\\
 B_r(x)&\sim\frac{\log_2x}{r(r-1)\log_3x}
       \quad(2\le r\le K).
\end{align*}
\end{theorem}

\begin{corollary}\label{cor:scalar}
The non-cyclic Sylow count satisfies
\[
 \frac{W-A_W(x)}{\sqrt{\frac{\log_2x}{2\log_3x}}}\Longrightarrow N(0,1),\qquad
 \Var_x W\sim\frac12\frac{\log_2x}{\log_3x}.
\]
The common-prime count $C(n)=\omega(\gcd(\varphi(n),\psi(n)))$ satisfies
\[
 \frac{C-A_C(x)}{\sqrt{\frac{6\log(9/8)\log_2x}{\log_3x}}}\Longrightarrow N(0,1),\qquad
 \Var_x C\sim6\log(9/8)\frac{\log_2x}{\log_3x}.
\]
In either limit one may use the exact arithmetic mean and
standard deviation.
\end{corollary}

Corollary \ref{cor:scalar} gives both the size of the fluctuations
and the proportions of integers for which they lie in a given range. If $\Phi$ denotes the standard normal distribution
function, then for every fixed pair of real numbers $a<b$,
\[
 \Pp_x\left(a<
       \frac{W(n)-A_W(x)}{\sqrt{\frac{\log_2x}{2\log_3x}}}\le b\right)
       \longrightarrow\Phi(b)-\Phi(a).
\]
In particular, the proportion of integers $n\le x$ whose count
lies within one or two asymptotic standard deviations of
$A_W(x)$ tends to approximately $68.27\%$ or $95.45\%$,
respectively. These limits describe the spread around the center, as well as
the typical size given by normal order.

\begin{corollary}[Normal orders]\label{cor:normalorders}
The normal orders of $W(n)$, $C(n)$ and $E_r(n)$, for each
fixed $r\ge2$, are respectively
\[
 \frac{\log_2n}{\log_3n},\qquad
 (2\log2)\frac{\log_2n}{\log_3n},\qquad
 \frac{\log_2n}{r(r-1)\log_3n}.
\]
For $W$ and $E_r$, these recover the results of
Pollack \cite{Pollack} and Morales--Polanco--Pollack \cite{MPP}.
\end{corollary}
\begin{proof}
Theorem \ref{thm:joint} gives the asserted mean asymptotics and
a variance $O_K(\frac{\log_2x}{\log_3x})$ for each count. Chebyshev's inequality
shows that each count divided by its mean tends to one in
probability. Since $\log_2n/\log_3n\sim \frac{\log_2x}{\log_3x}$ uniformly for
$\sqrt x\le n\le x$, and only $o(x)$ integers lie below
$\sqrt x$, the asserted normal orders follow.
\end{proof}

\begin{corollary}[Sum-of-divisors variant]\label{cor:sigma}
Replacing $C(n)$ by
$C_\sigma(n)=\omega(\gcd(\varphi(n),\sigma(n)))$, where
$\sigma(n)=\sum_{d\mid n}d$, leaves the
joint limit, centers and covariance constants in
Theorem \ref{thm:joint} unchanged.
\end{corollary}
\begin{proof}
Lemma \ref{lem:groups} gives $\|C_\sigma-C\|_{L^2(\Pp_x)}=O(1)$.
Slutsky's theorem transfers the joint limit. Cauchy--Schwarz,
with the $O_K(\sqrt{\frac{\log_2x}{\log_3x}})$ standard deviations from
Theorem \ref{thm:joint}, changes each covariance by
$O_K(\sqrt{\frac{\log_2x}{\log_3x}}+1)=o_K(\frac{\log_2x}{\log_3x})$, and the mean changes by $O(1)$.
\end{proof}

A first-order equivalent of a center need
not approximate it on the scale $\sqrt{\frac{\log_2x}{\log_3x}}$; we therefore retain
the explicit sums \eqref{eq:centers}.
For example,
\[
 \Var_x E_2\sim\frac7{16}\frac{\log_2x}{\log_3x},\quad
 \Var_x E_3\sim\frac7{48}\frac{\log_2x}{\log_3x},\quad
 \Cov_x(E_2,E_3)\sim-\frac1{32}\frac{\log_2x}{\log_3x}.
\]
The joint theorem also shows how these counts vary together.
The leading covariance between $E_r$ and $E_s$ is
negative for distinct fixed ranks, whereas that between $W$
and $E_r$ is positive. These signs describe how the counts fluctuate
around their means as the modulus varies. The normal orders of
the individual counts do not give this information.

\section{The model and arithmetic estimates}
Fix $K\ge2$. We first define the prime counts and their random
model. We then prove the estimates for truncating these counts
and comparing them with $W$, $C$, and $E_r$. These estimates
will be used in the proof of the main theorem in
Section~\ref{sec:proof}.

\subsection{Counts and the Bernoulli model}
For sufficiently large $x$, fix
\begin{equation}
 y=x^{1/\log_3x},\qquad Q=(\log_2x)^5,\qquad z=\frac{\log_2x}{(\log_3x)^2}.
 \label{eq:cutoffs}
\end{equation}
Here $y$ bounds input primes, $Q$ bounds output primes, and $z$
separates the almost deterministic small outputs.
These choices satisfy
\begin{equation}
 2<y<x,\qquad Q^4<y
 \label{eq:cutoffranges}
\end{equation}
for all sufficiently large $x$. Since $Q^4<y$, the tail estimates
start above the square of every modulus $q\ell\le Q^2$.

We count input prime divisors in the classes $1$ and $-1$
modulo $q$. This lets us treat all the statistics with the
same argument. For $q\ge5$, put
\[
 N_{q,\pm}(n)=\#\{p\mid n:p\equiv\pm1\pmod q\}.
\]
For $2\le k\le K+1$, the arithmetic occupancy counts are
\[
 \mathcal C_k(n)=\sum_{q\ge5}\ind_{\{N_{q,+}\ge k\}},\qquad
 \mathcal C_\cap(n)=\sum_{q\ge5}
       \ind_{\{N_{q,+}\ge1,N_{q,-}\ge1\}}.
\]
We index these types by $\mathcal H=\{2,3,\ldots,K+1,\cap\}$,
where $\cap$ denotes the two-sided condition.
Write $\mathcal C_h(n;y,Q)$ when inputs are restricted to
$p\le y$ and outputs to $5\le q\le Q$, and
$\mathcal C_h(n;\infty,Q)$ when only outputs are restricted.
Each sum is finite for each $n$.

In the random model, we replace each event $p\mid n$ by an
independent Bernoulli variable $B_p$, for $p\le y$, with
$\Pp(B_p=1)=1/p$. Their residue-class occupancies are
\[
 K_{q,\pm}=\sum_{\substack{p\le y\\p\equiv\pm1\pmod q}}B_p.
\]
For $2\le k\le K+1$, their threshold and two-sided indicators are
\[
 X_{k,q}=\ind_{\{K_{q,+}\ge k\}},\qquad
 X_{\cap,q}=\ind_{\{K_{q,+}\ge1,K_{q,-}\ge1\}}.
\]
Probability, expectation, variance, and covariance in this
model carry the subscript $B$.

For the cutoffs \eqref{eq:cutoffs}, let $n$ have law $\Pp_x$
and let $(G_p:p\le y)$ be independent geometric
variables with $\Pp(G_p=a)=(1-1/p)p^{-a}$ for
$a=0,1,\ldots$. Here $v_p(n)$ is the exponent of $p$ in $n$.
For total variation distance
$d_{\rm TV}(\mu,\nu)=\sup_A|\mu(A)-\nu(A)|$, put
\[
 \delta=\delta_x(y)
   =d_{\rm TV}\bigl((v_p(n):p\le y),(G_p:p\le y)\bigr).
\]
The integer Kubilius estimate of Tenenbaum
\cite[(1.7)]{Tenenbaum}, stated for the full vector of prime
valuations in \cite[Theorem A]{Ford}, applies since $2<y<x$
and $\log x/\log y=\log_3x$. For every fixed $\epsilon>0$, it gives
\begin{equation}
 \delta\ll_\epsilon (\log_3x)^{-\log_3x}+x^{-1+\epsilon}.
 \label{eq:TVinput}
\end{equation}
The total variation convention and uniform integer law in these
references agree with ours. Taking $\epsilon=1/2$, we use
\[
 (\log_3x)^{-\log_3x}
       =e^{-(\log_3x)(\log_4x)}
       =(\log_2x)^{-\log_4x}.
\]
It follows that
\begin{equation}
 (\log_2x)^A\delta=o(1)\quad\hbox{for every fixed }A>0.
 \label{eq:TVsmall}
\end{equation}
In particular, $Q^2\delta=o(\frac{\log_2x}{\log_3x})$ and $z^2\delta=o(1)$.
More explicitly, $Q^2=(\log_2x)^{10}$, so
\begin{align}
 \frac{\log_3x}{\log_2x}Q^2\delta
 &\ll(\log_3x)(\log_2x)^{9-\log_4x}\notag\\
 &\quad+(\log_3x)(\log_2x)^9x^{-1/2}\longrightarrow0.
 \label{eq:TVsecond}
\end{align}
Also $z^2\delta\le (\log_2x)^2\delta=o(1)$.
For each $G_p$, take the indicator that it is positive. These
indicators are independent and have mean $\Pp(G_p\ge1)=1/p$,
as in the model $B_p$.
A deterministic map of the entire input vector does not increase
total variation. Thus the whole joint vector
$(\mathcal C_h(n;y,Q):h\in\mathcal H)$ is approximated by
$\bigl(\sum_{5\le q\le Q}X_{h,q}:h\in\mathcal H\bigr)$
with total variation error at most $\delta$. There is no additional
factor for the number of output primes.
We will also use the corresponding expectation bound: for a
nonnegative bounded statistic, the difference of its expectations
is at most its maximum value times $\delta$. This follows by
integrating its tail probabilities under the two laws.

\subsection{Prime sums}
We next estimate the reciprocal-prime sums used to bound the
omitted inputs and the inputs shared by two output primes.
For $a\in\{1,-1\}$,
write
\[
 \Lambda_{d,a}(X)=\sum_{\substack{p\le X\\p\equiv a\pmod d}}\frac1p.
\]
\begin{lemma}\label{lem:AP}
Uniformly for every integer $d\ge2$ and $X\ge3$,
\begin{equation}
 \Lambda_{d,\pm1}(X)\ll
       \frac{1+\log\log(3X)}{\varphi(d)}.
 \label{eq:APall}
\end{equation}
Uniformly for primes $5\le q\le Q$,
\begin{equation}
 \Lambda_{q,\pm1}(y)=\frac{\log_2x-\log_4x}{q-1}
          +O\left(\frac{\log(2q)}q\right).
 \label{eq:APmain}
\end{equation}
For $d=q$ or $d=q\ell$, where $q,\ell\le Q$ are distinct
primes, and every reduced residue class $b\pmod d$,
\begin{equation}
 \sum_{\substack{y<p\le x\\p\equiv b\pmod d}}\frac1p
       \ll\frac{1+\log_4x}{\varphi(d)}.
 \label{eq:APtail}
\end{equation}
For distinct primes $q,\ell\in[z,Q]$ and every reduced class
$b\pmod{q\ell}$,
\begin{equation}
 \sum_{\substack{Q^4<p\le y\\p\equiv b\pmod{q\ell}}}\frac1p
       \ll\frac{\log_2x}{q\ell}.
 \label{eq:CRT}
\end{equation}
All constants are independent of the moduli and of $x$.
\end{lemma}
\begin{proof}
Write $\pi(t;d,b)$ for the number of primes $p\le t$ with
$p\equiv b\pmod d$. We use Brun--Titchmarsh in the form
\[
 \pi(t;d,b)\le\frac{2t}{\varphi(d)\log(t/d)}
       \quad(t>d,\ (b,d)=1).
\]
To prove \eqref{eq:APall} for the two special residues, we first
split the prime sum at $2d$.
In $2\le p\le\min(X,2d)$, the only candidates are $d-1,2d-1$
for $-1$, and $d+1$ for $+1$. We retain only primes not exceeding
$X$; in particular, the integer $1$ is excluded.
Their reciprocal sum is $O(1/d)$, which proves the estimate
when $X\le2d$. If $X>2d$, partial summation
on the remaining interval bounds it by
\[
 \frac{C_0}{\varphi(d)}\left\{\frac1{\log(X/d)}
             +\int_{2d}^X\frac{\dd t}{t\log(t/d)}\right\}
 \ll\frac{1+\log\log(3X)}{\varphi(d)}.
\]
This proves the all-modulus bound, using $\varphi(d)\le d$. The split at twice
the modulus also appears for the minus residue in
\cite[Lemma 2.4]{GuoRatio}. No all-modulus claim for arbitrary
reduced residues is being used.

For \eqref{eq:APmain}, apply Norton's uniform reciprocal-prime
estimate \cite[Lemma 1.22]{Norton} with upper limit $y$,
modulus $q$, and residue set $\{1\}$ or $\{q-1\}$.
Both sets consist of reduced residues because $q$ is prime.
The estimate holds for every $y\ge2$, uniformly in the modulus.
Our cutoffs give $y=x^{1/\log_3x}>Q^4>2$ for large $x$
and $5\le q\le Q=(\log_2x)^5$.
The correction for prime representatives below the modulus
vanishes, since $1$ is not prime and $q-1\ge4$ is even.
We therefore obtain
\[
 \Lambda_{q,\pm1}(y)=\frac{\log\log y}{q-1}
       +O\left(\frac{\log(3q)}{q-1}\right),
\]
with an absolute implied constant. Since
$\log\log y=\log_2x-\log_4x$, this gives \eqref{eq:APmain}.

For \eqref{eq:APtail}, $d\le Q^2$ and $y>Q^4\ge d^2$.
Thus Brun--Titchmarsh is $O(t/(\varphi(d)\log t))$ throughout
the interval. Partial summation, discarding the nonpositive
lower endpoint term, gives
\[
 \frac{C_0}{\varphi(d)}\left\{\frac1{\log x}
                      +\log\log x-\log\log y\right\}
       \ll\frac{1+\log_4x}{\varphi(d)}.
\]
For \eqref{eq:CRT} use the same argument on $(Q^4,y]$, where
$q\ell\le Q^2\le\sqrt t$, and
$\varphi(q\ell)\asymp q\ell$. See \cite{Davenport} for the
Brun--Titchmarsh inequality used here.
\end{proof}

\subsection{Removing large primes}
The model uses only input primes $p\le y$ and output primes
$q\le Q$. We next bound the error from the omitted primes,
both in mean and in mean square. For the output tail, we will
combine a small first moment with the following fourth-moment bound.

\begin{lemma}\label{lem:fourth}
For each fixed integer $m\ge1$,
\[
 \sum_{p\le x}\frac{[1+\omega(p-1)]^m}p\ll_m (\log_2x)^{m+1},
 \qquad \E_x\omega(\varphi(n))^4\ll (\log_2x)^8.
\]
\end{lemma}
\begin{proof}
For a fixed integer $j\ge1$, choosing $j$ distinct prime factors
of $p-1$ is the same as choosing a squarefree divisor with $j$
prime factors. This expansion and \eqref{eq:APall} give
\begin{align*}
 \sum_{p\le x}\frac1p\binom{\omega(p-1)}j
 &=\sum_{\substack{d\le x-1\\d\ \text{squarefree}\\\omega(d)=j}}
       \qquad\sum_{\substack{p\le x\\p\equiv1\pmod d}}\frac1p\\
 &\ll \log_2x\sum_{\substack{d\le x-1\\d\ \text{squarefree}\\\omega(d)=j}}
                  \frac1{\varphi(d)}
 \le\frac{\log_2x}{j!}\left(\sum_{q\le x}\frac1{q-1}\right)^j
 \ll_j (\log_2x)^{j+1}.
\end{align*}
For $j=0$ the bound is $O(\log_2x)$. Since
$\binom{t}{0},\ldots,\binom{t}{m}$ form a basis for polynomials
of degree at most $m$, there are constants $c_{m,j}$ such that
\[
 (1+t)^m=\sum_{j=0}^m c_{m,j}\binom tj.
\]
The preceding bounds therefore give
$\sum_{p\le x}[1+\omega(p-1)]^m/p
 \ll_m\sum_{j=0}^m (\log_2x)^{j+1}\ll_m (\log_2x)^{m+1}$.
Put $b_p=1+\omega(p-1)$ and $F(n)=\sum_{p\mid n}b_p$.
Then $\omega(\varphi(n))\le F(n)$. Write
$D_j=\sum_{p\le x}b_p^j/p\ll (\log_2x)^{j+1}$, for $1\le j\le4$.
For distinct input primes $p_1,\ldots,p_t$ we have
$\Pp_x(p_1\cdots p_t\mid n)\le1/(p_1\cdots p_t)$.
Expand the fourth power and group the input indices according
to a set partition $\pi$ of $\{1,2,3,4\}$. Each block of $\pi$
corresponds to one input prime. Apply the preceding divisibility
bound to the distinct primes assigned to the blocks, and then
drop the distinctness restriction in the nonnegative sums. This gives
\[
 \E_xF^4
 \ll D_1^4+6D_2D_1^2+3D_2^2+4D_3D_1+D_4
 \ll (\log_2x)^8.
\]
The five terms correspond, in order, to block sizes
$1+1+1+1$, $2+1+1$, $2+2$, $3+1$, and $4$.
The corresponding multiplicities are $1,6,3,4,1$, respectively.
No independence of arithmetic divisibility events is used.
\end{proof}
We now bound the two tails.

\begin{lemma}\label{lem:trunc}
Simultaneously for $h\in\mathcal H$,
\begin{align}
 \E_x|\mathcal C_h-\mathcal C_h(n;y,Q)|
       &\ll_K(1+\log_4x)^2,\label{eq:Lone}\\
 \E_x|\mathcal C_h-\mathcal C_h(n;y,Q)|^2
       &\ll_K(1+\log_4x)^4
          +\left(\frac{\log_2x}{\log_3x}\right)^{2/3}\notag\\
       &=o\left(\frac{\log_2x}{\log_3x}\right).
       \label{eq:Ltwo}
\end{align}
\end{lemma}
\begin{proof}
\emph{Output primes.}
Every nonzero output indicator needs at least two distinct input
primes. In the cap case their sign classes are disjoint.
For distinct inputs $p,t$,
$\Pp_x(pt\mid n)\ll1/(pt)$. Summing over eligible pairs
bounds each output probability by
$\Lambda_{q,+1}(x)^2$ in the threshold case and by
$\Lambda_{q,+1}(x)\Lambda_{q,-1}(x)$ in the cap case.
Both are $O((\log_2x)^2/q^2)$ by \eqref{eq:APall}.
Thus the output tail
$R_h=\mathcal C_h-\mathcal C_h(n;\infty,Q)\ge0$ satisfies
\begin{align}
 \E_xR_h&\ll(\log_2x)^2\sum_{q>Q}q^{-2}\notag\\
 &\ll\frac{(\log_2x)^2}{Q\log Q}
       \ll\frac1{(\log_2x)^3\log_3x}.
 \label{eq:outmean}
\end{align}
This uses the all-modulus bound, including arbitrarily large
output primes. Each such output divides $\varphi(n)$, so
$R_h\le\omega(\varphi(n))$. Apply H\"older with exponents
$3/2$ and $3$ to $R_h^{2/3}R_h^{4/3}$. Lemma \ref{lem:fourth}
then gives
\begin{align}
 \E_xR_h^2
 &\le(\E_xR_h)^{2/3}(\E_xR_h^4)^{1/3}\notag\\
 &\ll\left(\frac1{(\log_2x)^3\log_3x}\right)^{2/3}
            \bigl((\log_2x)^8\bigr)^{1/3}\notag\\
 &=\left(\frac{\log_2x}{\log_3x}\right)^{2/3}
       =o\left(\frac{\log_2x}{\log_3x}\right).
 \label{eq:outsecond}
\end{align}

\emph{Input primes.}
For the input tail and $q\le \log_2x$, a changed indicator requires
at least one tail witness, of expected number $O((1+\log_4x)/q)$ by
\eqref{eq:APtail}. For $\log_2x<q\le Q$, a changed indicator requires
a tail witness $p>y$ and a distinct eligible witness $t\le x$.
Summing over these pairs, using
$\Pp_x(pt\mid n)\le1/(pt)$ together with
\eqref{eq:APtail} and \eqref{eq:APall}, bounds its expected change
by $O_K((1+\log_4x)\log_2x/q^2)$.
Consequently the total expected input error is at most
\[
 C_K(1+\log_4x)\left\{\sum_{q\le \log_2x}\frac1q
                    +\log_2x\sum_{q>\log_2x}\frac1{q^2}\right\}
       \ll_K(1+\log_4x)^2.
\]
For its second moment put
\[
 \kappa_Q(p)=\#\{5\le q\le Q:q\mid(p-1)(p+1)\},\qquad
 J(n)=\sum_{\substack{p>y\\p\mid n}}\kappa_Q(p).
\]
Every changed output has an eligible tail witness. Therefore
$0\le\mathcal C_h(n;\infty,Q)-\mathcal C_h(n;y,Q)\le J$ for every $h$.
The two sign classes are disjoint for $q\ge5$. In the square
of $\kappa_Q(p)$, a diagonal term therefore uses a single
modulus $q$. For distinct $q,\ell$, the eligible inputs lie in
at most four reduced classes modulo $q\ell$, by the Chinese
remainder theorem (CRT). Applying \eqref{eq:APtail} gives
\begin{align*}
 \sum_{y<p\le x}\frac{\kappa_Q(p)}p
 &\ll(1+\log_4x)\sum_{5\le q\le Q}\frac1q
       \ll(1+\log_4x)^2,\\
 \sum_{y<p\le x}\frac{\kappa_Q(p)^2}p
 &\ll(1+\log_4x)\left\{\sum_{5\le q\le Q}\frac1q
       +\left(\sum_{5\le q\le Q}\frac1q\right)^2\right\}\\
 &\ll(1+\log_4x)^3.
\end{align*}
To separate equal from distinct input primes, expand directly:
\[
 J^2=\sum_{y<p\le x}\kappa_Q(p)^2\ind_{\{p\mid n\}}
       +2\sum_{y<p<p'\le x}
          \kappa_Q(p)\kappa_Q(p')\ind_{\{pp'\mid n\}}.
\]
Applying the divisibility bounds to these two sums gives
\begin{equation}
 \E_xJ^2
 \ll\sum_{y<p\le x}\frac{\kappa_Q(p)^2}{p}
       +\left(\sum_{y<p\le x}\frac{\kappa_Q(p)}p\right)^2
 \ll(1+\log_4x)^4.
 \label{eq:inputsecond}
\end{equation}
Finally,
\[
 \mathcal C_h-\mathcal C_h(n;y,Q)
 =R_h+\bigl(\mathcal C_h(n;\infty,Q)-\mathcal C_h(n;y,Q)\bigr).
\]
The first-moment bound follows by addition. For the second
moment, use $(a+b)^2\le2a^2+2b^2$ with
\eqref{eq:outsecond} and \eqref{eq:inputsecond}.
\end{proof}

\subsection{Comparison with the group counts}
The truncation bounds concern counts in residue classes.
We now compare these counts with $W$, $C$, and $E_r$.
The identity
\[
 \ind_{\{N_{q,+}=r\}}=
       \ind_{\{N_{q,+}\ge r\}}-\ind_{\{N_{q,+}\ge r+1\}}
\]
shows that $\mathcal C_r-\mathcal C_{r+1}$ counts output primes
with exactly $r$ witnesses. To compare this count with $E_r$,
we use the local group structure.
For an odd output prime $q$, each local Sylow $q$-factor is
cyclic. Write $p^a\parallel n$ when $v_p(n)=a$.
For $p^a\parallel n$ with $p\ne q$, this factor is
nontrivial exactly when $q\mid p-1$; for $p=q$, it is
nontrivial exactly when $a\ge2$.
The rank of the full Sylow $q$-subgroup, meaning its
minimal number of generators, is
\begin{equation}
 N_{q,+}(n)+\ind_{\{q^2\mid n\}}.
 \label{eq:rankidentity}
\end{equation}
This rank also equals the number of nontrivial cyclic
$q$-factors in the local decomposition.
The Sylow subgroup is non-elementary precisely when $q^3\mid n$ or when
$q^2\mid p-1$ for some prime $p\mid n$.
These statements follow from the standard local decomposition
of $(\mathbb Z/n\mathbb Z)^\times$; see also \cite[Lemma 2.1]{MPP}.

\begin{lemma}\label{lem:groups}
The errors $W-\mathcal C_2$ and $C-\mathcal C_\cap$ have
bounded $L^2(\Pp_x)$ norm. For $2\le r\le K$,
\begin{equation}
 \E_x|E_r-(\mathcal C_r-\mathcal C_{r+1})|^2
       \ll_K (\log_3x)^2.
 \label{eq:elementaryLtwo}
\end{equation}
The difference between $C(n)$ and
$\omega(\gcd(\varphi(n),\sigma(n)))$ also has bounded $L^2$ norm.
\end{lemma}
\begin{proof}
\emph{The counts $W$ and $C$.}
For $q\ge5$, \eqref{eq:rankidentity} shows that the indicator
in $W$ differs from $\ind_{\{N_{q,+}\ge2\}}$ precisely when
$N_{q,+}=1$ and $q^2\mid n$: the full rank is then two
although there is only one progression witness. Hence
\[
 0\le W-\mathcal C_2\le2+\sum_{q\ge5}\ind_{\{q^2\mid n\}};
\]
the output primes $2,3$ contribute at most two. Put
$\rad(n)=\prod_{p\mid n}p$. The identities
\[
 \varphi(n)=\frac{n}{\rad(n)}\prod_{p\mid n}(p-1),\qquad
 \psi(n)=\frac{n}{\rad(n)}\prod_{p\mid n}(p+1)
\]
show, for $q\ge5$, that $q$ divides their gcd precisely when
$q^2\mid n$ or both sign classes have a witness. Thus
$0\le C-\mathcal C_\cap\le2+\sum_{q\ge5}\ind_{\{q^2\mid n\}}$
as well. Since
\[
 \E_x\left(\sum_{q\ge5}\ind_{\{q^2\mid n\}}\right)^2
       \ll\sum_q q^{-2}+\sum_{q\ne\ell}q^{-2}\ell^{-2}\ll1,
\]
both differences are bounded in $L^2$.

\emph{Fixed-rank elementary subgroups.}
First consider $5\le q<z$.
If either the elementary-rank indicator or the exact-occupancy
indicator is nonzero, then $N_{q,+}\le r\le K$.
Restricting inputs to $p\le y$ only enlarges this event.
Uniformly for $q<z$, \eqref{eq:APmain} gives
$\Lambda_{q,+1}(y)\gg \log_2x/q\gg (\log_3x)^2$.
In the following products and sums, all input primes satisfy
$p\le y$ and $p\equiv1\pmod q$.
Expansion according to the set of successful inputs gives,
for every $j\ge0$,
\begin{align*}
 \Pp_B(K_{q,+}=j)
 &=\prod_p(1-1/p)
       \sum_{p_1<\cdots<p_j}\prod_{i=1}^j\frac1{p_i-1}\\
 &\le\frac{e^{-\Lambda_{q,+1}(y)}}{j!}
             \left(\sum_p\frac1{p-1}\right)^j\\
 &\le e^{-\Lambda_{q,+1}(y)}
             \frac{(2\Lambda_{q,+1}(y))^j}{j!}.
\end{align*}
The empty product is one. We used $1-t\le e^{-t}$,
$1/(p-1)\le2/p$, and the elementary symmetric-sum bound
obtained by expanding the $j$th power of a nonnegative sum.
Summing over $0\le j\le K$ yields
\[
 \Pp_B(K_{q,+}\le K)
 \le e^{-\Lambda_{q,+1}(y)}
          \sum_{j=0}^K\frac{(2\Lambda_{q,+1}(y))^j}{j!}
 \ll_K e^{-c_K(\log_3x)^2}.
\]
Here $c_K>0$ depends only on $K$.
Write $N^y_{q,+}$ for the count restricted to inputs $p\le y$,
and put $I_q^y(n)=\ind_{\{N^y_{q,+}(n)\le K\}}$.
The small-output discrepancy is at most
$\sum_{5\le q<z}I_q^y$. Its square depends only on inputs
$p\le y$ and is bounded by $z^2$. We can therefore apply
\eqref{eq:TVinput} to this squared sum. In the Bernoulli model,
the square is at most $z^2$ times the indicator that at least
one low-occupancy event occurs. The union bound then gives
\[
 O_K\{z^3 e^{-c_K(\log_3x)^2}+z^2\delta\}=o(1)
\]
for the small-output contribution to \eqref{eq:elementaryLtwo}.

For $q\ge z$, the elementary and exact-occupancy indicators
agree unless $q^2\mid n$ or $q^2\mid p-1$ for some $p\mid n$.
Let
\[
 \kappa_z(p)=\#\{q\ge z:q^2\mid p-1\},\qquad
 J_z(n)=\sum_{p\mid n}\kappa_z(p).
\]
The large-output discrepancy is at most
$\sum_{q\ge5}\ind_{\{q^2\mid n\}}+J_z$.
Using \eqref{eq:APall} with $q^2$ and $q^2\ell^2$, respectively,
we obtain
\begin{align*}
 \sum_{p\le x}\frac{\kappa_z(p)}p&
       \ll \log_2x\sum_{q\ge z}q^{-2}\ll\frac{\log_2x}{z\log z}\ll \log_3x,\\
 \sum_{p\le x}\frac{\kappa_z(p)^2}p&
       \ll \log_2x\left\{\sum_{q\ge z}q^{-2}
                  +\left(\sum_{q\ge z}q^{-2}\right)^2\right\}
       \ll \log_3x.
\end{align*}
The mixed term uses $\varphi(q^2\ell^2)\asymp q^2\ell^2$ for
distinct output primes. For input primes, the same expansion
as in \eqref{eq:inputsecond} gives
\begin{align*}
 \E_xJ_z^2
 &\ll\sum_{p\le x}\frac{\kappa_z(p)^2}{p}
       +\sum_{\substack{p,t\le x\\p\ne t}}
                    \frac{\kappa_z(p)\kappa_z(t)}{pt}\\
 &\le\sum_{p\le x}\frac{\kappa_z(p)^2}{p}
       +\left(\sum_{p\le x}\frac{\kappa_z(p)}p\right)^2\\
 &\ll (\log_3x)^2.
\end{align*}
Together with the small-output bound and the omitted primes
$2,3$, this proves \eqref{eq:elementaryLtwo}.

\emph{The sum-of-divisors variant.}
Finally, set
\[
 \Delta(n)=\sum_q|v_q(\sigma(n))-v_q(\psi(n))|\log q.
\]
The local factors agree when the input prime has exponent one, since
$\sigma(p)=\psi(p)=p+1$. For $a\ge2$,
\[
 \sigma(p^a)=\frac{p^{a+1}-1}{p-1}\le2p^a,\qquad
 \psi(p^a)=p^{a-1}(p+1)\le2p^a.
\]
Consequently,
\[
 \sum_q|v_q(\sigma(p^a))-v_q(\psi(p^a))|\log q
 \le\log\sigma(p^a)+\log\psi(p^a)\ll a\log p.
\]
Multiplicativity and the triangle inequality for valuations give
\[
 \Delta(n)\ll\sum_{\substack{p^a\parallel n\\a\ge2}}a\log p.
\]
It suffices to bound the second moment. For a fixed input
prime, events with different exact exponents are disjoint.
For distinct input primes $p,t$, product divisibility bounds
$\Pp_x(p^at^b\mid n)$ by $p^{-a}t^{-b}$. Expanding the
square therefore bounds the second moment by a constant times
\[
 \sum_p\sum_{a\ge2}\frac{a^2(\log p)^2}{p^a}
       +\left(\sum_p\sum_{a\ge2}\frac{a\log p}{p^a}\right)^2<\infty.
\]
Here $\sum_{a\ge2}a^2p^{-a}\ll p^{-2}$ uniformly for $p\ge2$,
and both $\sum_p(\log p)^2/p^2$ and $\sum_p(\log p)/p^2$
converge. The first moment is bounded by Cauchy--Schwarz.
A prime can divide one gcd but not the other only if it divides
one of the two second arguments but not the other. Thus the two
gcd prime counts differ by at most $\Delta(n)/\log2$.
This proves the last assertion.
\end{proof}

\section{Proof of Theorem \ref{thm:joint}}\label{sec:proof}
We use the model and arithmetic estimates from the preceding
section, with $K$ fixed. We first find the means and control the
dependence between different output primes. We then compute
the covariance matrix and prove the Gaussian limit in the model.
Finally, we return to the arithmetic counts using the bounds in
Lemmas~\ref{lem:trunc} and~\ref{lem:groups}.

\subsection{Means and dependence in the model}
We begin with the means in the Bernoulli model. For $a\ge0$
and integers $k,r\ge2$, put
\[
 H_k(a)=1-e^{-a}\sum_{j=0}^{k-1}\frac{a^j}{j!},\qquad
 H_\cap(a)=(1-e^{-a})^2,\qquad
 P_r(a)=e^{-a}\frac{a^r}{r!}.
\]
These are the Poisson probabilities of at least $k$ witnesses,
of a witness in each of two independent sign classes, and of
exactly $r$ witnesses, respectively. For $h\in\mathcal H$,
$H_h$ denotes the corresponding threshold or two-sided probability.
We will use the elementary identity
\begin{equation}
 H_r(a)-H_{r+1}(a)=P_r(a).
 \label{eq:centerdifference}
\end{equation}

\emph{Means.}
We first compute the mean of each model count.
Le Cam's inequality \cite[Proposition 1]{LeCam} gives, for each
sign class, total variation error at most
\[
 \sum_{\substack{p\le y\\p\equiv\pm1\pmod q}}p^{-2}\ll q^{-2}
\]
when $K_{q,\pm}$ is replaced by a Poisson variable of mean
$\Lambda_{q,\pm1}(y)$.
The sign classes are disjoint for $q\ge5$, so their model counts
are independent. In the cap case the two Poisson means need
not coincide; for $s,t,a\ge0$ we use
\[
 \left|(1-e^{-s})(1-e^{-t})-(1-e^{-a})^2\right|
       \le |s-a|+|t-a|.
\]
This follows from $0\le1-e^{-v}\le1$ and
$|e^{-v}-e^{-w}|\le|v-w|$ for $v,w\ge0$.
Together with the bounded derivatives of the fixed threshold
probabilities and \eqref{eq:APmain}, this gives
\[
 |\E_BX_{h,q}-H_h(\log_2x/(q-1))|
       \ll_K\frac{\log_4x+\log(2q)}q+q^{-2}.
\]
Summing over $q\le Q$ gives $O_K(\log_3x)$, since
$\sum_{q\le Q}(\log q)/q\ll\log Q$ and
$(\log_4x)\sum_{q\le Q}1/q\ll(\log_4x)^2=o(\log_3x)$.
The center tail is $O_K((\log_2x)^2/(Q\log Q))=o(1)$, as
$H_h(a)\ll_K a^2$ at zero. Thus, for every $h\in\mathcal H$,
\begin{equation}
 \E_B\sum_{5\le q\le Q}X_{h,q}
   =\sum_{q\ge5}H_h\left(\frac{\log_2x}{q-1}\right)+O_K(\log_3x).
 \label{eq:modelmean}
\end{equation}

\emph{Small output primes.}
We next show that the small output primes contribute almost
constant values.
For $q<z$, the symmetric-sum estimate in Lemma \ref{lem:groups}
bounds the failure probability of any fixed threshold by
$O_K(e^{-c_K(\log_3x)^2})$. For the cap coordinate its failure
probability is at most
$e^{-\Lambda_{q,+1}(y)}+e^{-\Lambda_{q,-1}(y)}$.
A union bound and a deterministic square bound yield
$O_K(z^3e^{-c_K(\log_3x)^2})=o(1)$. Hence, for every
$h\in\mathcal H$,
\begin{equation}
 \E_B\left|\sum_{5\le q<z}X_{h,q}
        -\#\{q:5\le q<z\}\right|^2=o_K(1).
 \label{eq:smalloutputs}
\end{equation}

\emph{Different output primes.}
We next estimate the covariances between different output primes.
These counts can depend on the same input prime. To obtain a bound, write
$B=(B_1,\ldots,B_m)$ for a finite vector of independent
Bernoulli variables, with $\theta_i=\Pp(B_i=1)$.
The notation $B_{\ne i}$ means the vector with coordinate $i$
removed. For increasing functions $F,G$ of $B$, write
$D_iF=F|_{B_i=1}-F|_{B_i=0}$. Then
\begin{equation}
 0\le\Cov(F,G)\le\sum_i\theta_i
             \sqrt{\E[(D_iF)^2]\,\E[(D_iG)^2]}.
 \label{eq:covbound}
\end{equation}
Indeed, let $B'$ be an independent copy of $B$, and replace
the coordinates of $B'$ successively by those of $B$ in the
argument of $G$. At step $i$, put
\[
 C_i=(B_1,\ldots,B_{i-1},B'_{i+1},\ldots,B'_m).
\]
Conditioning on $B_{\ne i}$ and $B'_{\ne i}$ gives
\begin{align*}
 &\E\bigl[F(B)G(B_1,\ldots,B_i,B'_{i+1},\ldots,B'_m)\bigr]\\*
 &\quad-\E\bigl[F(B)G(B_1,\ldots,B_{i-1},B'_i,\ldots,B'_m)\bigr]\\
 &=\theta_i(1-\theta_i)
       \E\bigl[D_iF(B_{\ne i})D_iG(C_i)\bigr].
\end{align*}
Each argument has the same product distribution as the other
coordinates of the input vector. Summing the increments and
applying Cauchy--Schwarz gives \eqref{eq:covbound}.
For increasing $F,G$, both discrete derivatives are nonnegative.
The identity therefore also proves association for increasing
functions of the input vector, including the threshold indicators.

Fix $h,j\in\mathcal H$. We now apply this bound to the
input variables $B_p$.
Here $D_p$ is the same difference operator for coordinate $B_p$.
A nonzero derivative of $X_{k,q}$ requires exactly $k-1$
other witnesses, and $k\ge2$; a cap derivative requires another
witness in the other sign class. These derivatives take values
zero or one. Their second moments are therefore bounded by
the probability of at least one such other witness, which is
at most the expected number of those witnesses. By
\eqref{eq:APall}, that expectation is $O_K(\log_2x/q)$.
Thus, for every type,
\begin{equation}
 \E_B[(D_pX_{h,q})^2]\le\min(1,C_K \log_2x/q)=:w_q.
 \label{eq:pivotal}
\end{equation}
The derivative is zero unless the input prime belongs to a
sign class used by the indicator. We do not need the derivatives
themselves to be increasing.
We split the input primes at $Q^4$. For $p\le Q^4$, put
\[
 \kappa(p)=\#\{z\le q\le Q:q\mid(p-1)(p+1)\}.
\]
This counts all output primes that this input can affect.
It satisfies
\[
 \kappa(p)\le\frac{\log((p-1)(p+1))}{\log z}
       \le\frac{8\log Q}{\log z}=O(1).
\]
An input with $\kappa(p)>0$ satisfies $p\ge z-1$.
Sum \eqref{eq:covbound} over distinct outputs $q,\ell$ and
interchange the finite nonnegative sums. A fixed input $p$
can contribute to at most $\kappa(p)^2$ ordered output pairs.
Since $w_q\le1$, the part from $p\le Q^4$ is bounded by
\begin{align*}
 \sum_{p\le Q^4}\frac{\kappa(p)^2}{p}
 &\ll_K\sum_{z-1\le p\le Q^4}\frac1p\\
 &\ll_K 1+\log\frac{\log(Q^4)}{\log(z-1)}=O_K(1).
\end{align*}
The last bound uses $\log(Q^4)/\log(z-1)\to20$.
This grouping handles mixed signs even when the first prime
in a CRT class is small compared with its modulus.

For $Q^4<p\le y$ and fixed distinct $q,\ell\in[z,Q]$, a
shared input must satisfy
\[
 p\equiv\pm1\pmod q,\qquad p\equiv\pm1\pmod\ell.
\]
The four choices of signs give four reduced classes modulo
$q\ell$, by the Chinese remainder theorem. Applying
\eqref{eq:CRT} to each class gives
\[
 \sum_{\substack{Q^4<p\le y\\q\mid(p-1)(p+1)\\
                    \ell\mid(p-1)(p+1)}}\frac1p
       \ll\frac{\log_2x}{q\ell}.
\]
Together with \eqref{eq:covbound} and \eqref{eq:pivotal}, this
bounds the part from these inputs by
\begin{align*}
 &C_K\log_2x\sum_{\substack{z\le q,\ell\le Q\\q\ne\ell}}
                 \frac{\sqrt{w_qw_\ell}}{q\ell}\\
 &\qquad\le C_K\log_2x\left(\sum_{z\le q\le Q}
                              \frac{\sqrt{w_q}}q\right)^2.
\end{align*}
Here we added the nonnegative diagonal terms after applying
the CRT estimate to distinct $q,\ell$. To bound the last sum,
split at $C_K\log_2x$, using the constant in \eqref{eq:pivotal}:
\begin{align*}
 \sum_{z\le q\le C_K\log_2x}\frac{\sqrt{w_q}}q
 &\le\sum_{z\le q\le C_K\log_2x}\frac1q
       \ll_K\frac{\log_4x}{\log_3x},\\
 \sum_{C_K\log_2x<q\le Q}\frac{\sqrt{w_q}}q
 &\le\sqrt{C_K\log_2x}\sum_{q>C_K\log_2x}q^{-3/2}
       \ll_K\frac1{\log_3x}.
\end{align*}
Combining the two ranges gives, for every $h,j\in\mathcal H$,
\begin{align}
 &\sum_{\substack{z\le q,\ell\le Q\\q\ne\ell}}
                  \Cov_B(X_{h,q},X_{j,\ell})\notag\\
 &\qquad\ll_K1+\frac{(\log_2x)(\log_4x)^2}{(\log_3x)^2}
       =o_K\left(\frac{\log_2x}{\log_3x}\right).
 \label{eq:mixedcov}
\end{align}
Every summand is nonnegative by association. Thus the total
covariance between distinct output primes is $o_K(\frac{\log_2x}{\log_3x})$.

\subsection{The model covariance matrix}
We now compute the covariance contributed by each output
prime and sum over $q$. The model coordinates corresponding
to $W$, $C$, and $E_r$ are
\begin{align*}
 Z_{0,q}&=X_{2,q},\qquad Z_{1,q}=X_{\cap,q},\\
 Z_{r,q}&=X_{r,q}-X_{r+1,q}=\ind_{\{K_{q,+}=r\}}
       \quad(2\le r\le K).
\end{align*}
For $0\le i\le K$, put $Z_i=\sum_{5\le q\le Q}Z_{i,q}$.
Write $Z=(Z_0,\ldots,Z_K)$.
By \eqref{eq:modelmean} and \eqref{eq:centerdifference},
\begin{align}
 \E_BZ_0&=A_W(x)+O_K(\log_3x),\notag\\
 \E_BZ_1&=A_C(x)+O_K(\log_3x),\notag\\
 \E_BZ_r&=B_r(x)+O_K(\log_3x)\qquad(2\le r\le K).
 \label{eq:Zmean}
\end{align}
Equation \eqref{eq:smalloutputs} gives a constant approximation
to each small-output sum. This constant is the number of small
output primes in coordinates $0,1$, and zero in the exact-rank
coordinates. The $L^2$ error is $o(1)$, so each centered
small-output sum has $L^2$ norm $o(1)$.

\emph{Covariance limit.}
For $q$ of size $\log_2x$, the two residue-class counts are
approximated by independent Poisson variables. Let $U,V$ be
independent Poisson variables of common mean $a>0$, and define
\[
 f_0=\ind_{\{U\ge2\}},\quad f_1=\ind_{\{U\ge1,V\ge1\}},\quad
 f_r=\ind_{\{U=r\}}\quad(2\le r\le K).
\]
Covariances involving the $f_i$ are taken under this Poisson law.
Fix $0<\alpha<\beta<\infty$. On the interval
$\alpha\log_2x\le q\le\beta\log_2x$, the model
block law converges uniformly to this law with $a=\log_2x/q$ by
\eqref{eq:APmain} and the Poisson approximation. There are
$O_{\alpha,\beta}(\log_2x/\log_3x)$ primes in this interval.
Thus the uniform $o(1)$ error in each covariance is still
$o(1)$ after summing and multiplying by $\log_3x/\log_2x$.
For fixed $0\le i,j\le K$, write $g(a)=\Cov(f_i,f_j)$.
Partial summation with the prime number theorem gives
\begin{align*}
 &\frac{\log_3x}{\log_2x}
       \sum_{\alpha \log_2x\le q\le\beta \log_2x}
       \Cov_B(Z_{i,q},Z_{j,q})\\
 &\quad=\frac{\log_3x}{\log_2x}
       \int_{\alpha \log_2x}^{\beta \log_2x}g(\log_2x/t)\frac{\dd t}{\log t}
       +o_{K,\alpha,\beta}(1)\\
 &\quad=\int_{1/\beta}^{1/\alpha}
       g(a)\frac{\log_3x}{\log_3x-\log a}\frac{\dd a}{a^2}
       +o_{K,\alpha,\beta}(1).
\end{align*}
The substitution is $a=\log_2x/t$, and the logarithmic ratio tends
uniformly to one on this fixed interval. Bounds for small and
large $q$, given below, then extend the limit to
\begin{equation}
 \frac{\log_3x}{\log_2x}\sum_{5\le q\le Q}\Cov_B(Z_{i,q},Z_{j,q})
       \longrightarrow
       \int_0^\infty\Cov(f_i,f_j)\frac{\dd a}{a^2}.
 \label{eq:integral}
\end{equation}
For fixed $\alpha>0$, the number
of primes below $\alpha \log_2x$, divided by $\frac{\log_2x}{\log_3x}$, has limsup
$O(\alpha)$, and each covariance is bounded. For $q>\beta \log_2x$
every nonzero coordinate requires two witnesses, so its
probability is $O_K((\log_2x)^2/q^2)$ and the tail divided by $\frac{\log_2x}{\log_3x}$ is
$O_K(1/\beta)$. We first let $x\to\infty$, and then let
$\alpha\downarrow0$, $\beta\uparrow\infty$. This proves
\eqref{eq:integral}, also when the covariance is negative.

\emph{Computing the constants.}
We now evaluate the integrals in \eqref{eq:integral}.
For $2\le r,s\le K$, the local covariances are
\begin{align*}
 \Cov(f_r,f_s)&=\ind_{\{r=s\}}P_r(a)-P_r(a)P_s(a),\\
 \Cov(f_0,f_r)&=P_r(a)(1+a)e^{-a},\\
 \Cov(f_1,f_r)&=P_r(a)(e^{-a}-e^{-2a}),\\
 \Cov(f_0,f_1)&=H_2(a)e^{-a}(1-e^{-a}).
\end{align*}
To evaluate the first three, use
\[
 \int_0^\infty a^m e^{-ca}\dd a=\frac{m!}{c^{m+1}}
       \qquad(m=0,1,\ldots,\ c>0).
\]
This identity follows by repeated integration by parts. It gives
\begin{align*}
 \int_0^\infty\frac{P_r(a)}{a^2}\dd a
       &=\frac{(r-2)!}{r!}=\frac1{r(r-1)},\\
 \int_0^\infty\frac{P_r(a)P_s(a)}{a^2}\dd a
       &=\frac{(r+s-2)!}{2^{r+s-1}r!s!}.
\end{align*}
Substituting these integrals into the first covariance formula
gives the entries for two fixed-rank counts.
Next, we use the same identity to compute the mixed entries:
\begin{align*}
 \int_0^\infty\frac{\Cov(f_0,f_r)}{a^2}\dd a
       &=\frac{(r-2)!}{2^{r-1}r!}
         +\frac{(r-1)!}{2^rr!}
         =\frac{r+1}{2^r r(r-1)},\\
 \int_0^\infty\frac{\Cov(f_1,f_r)}{a^2}\dd a
       &=\frac{(r-2)!}{r!}\bigl(2^{1-r}-3^{1-r}\bigr).
\end{align*}
These are the corresponding entries in \eqref{eq:covmatrix}.

For finite real coefficients $c_\nu$ and positive numbers
$b_\nu$, Frullani's identity and integration by parts give:
\begin{align*}
 \int_0^\infty\frac{\sum_\nu c_\nu e^{-b_\nu a}}a\dd a
   &=-\sum_\nu c_\nu\log b_\nu
       &&\left(\sum_\nu c_\nu=0\right),\\
 \int_0^\infty\frac{\sum_\nu c_\nu e^{-b_\nu a}}{a^2}\dd a
   &=\sum_\nu c_\nu b_\nu\log b_\nu
       &&\left(\sum_\nu c_\nu=\sum_\nu c_\nu b_\nu=0\right).
\end{align*}
The two zero-sum conditions in the second line make the numerator
$O(a^2)$ at zero, so the boundary terms vanish.
Applying these formulas to the last local covariance gives
\begin{align*}
 \int_0^\infty\frac{H_2(a)e^{-a}(1-e^{-a})}{a^2}\dd a
 &=\int_0^\infty\frac{e^{-a}-2e^{-2a}+e^{-3a}}{a^2}\dd a
       +\int_0^\infty\frac{-e^{-2a}+e^{-3a}}a\dd a\\
 &=(3\log3-4\log2)+\log(2/3)=\log(9/8).
\end{align*}
For the two remaining diagonal entries, note that
$H_2(a)=a^2/2+O(a^3)$ at zero, $H_2'(a)=ae^{-a}$, and
$H_2$ is bounded at infinity. Thus the boundary terms in the
following integrations by parts also vanish:
\[
 \int_0^\infty\frac{H_2(a)}{a^2}\dd a=1,\qquad
 \int_0^\infty\frac{H_2(a)^2}{a^2}\dd a
       =2\int_0^\infty H_2(a)e^{-a}\dd a=\frac12,
\]
and
\[
 H_\cap(a)(1-H_\cap(a))
       =2e^{-a}-5e^{-2a}+4e^{-3a}-e^{-4a}.
\]
The second exponential-sum identity above therefore gives
\[
 \int_0^\infty\frac{H_\cap(a)(1-H_\cap(a))}{a^2}\dd a
       =-10\log2+12\log3-4\log4=6\log(9/8).
\]
Subtracting the two $H_2$ integrals gives the first diagonal entry.
We have now obtained all entries in \eqref{eq:covmatrix}.

The same prime-sum argument applies to the means, since the
local probabilities are bounded and are $O_K(a^2)$ at zero.
For the two-sided count, integration by parts gives
\[
 \int_0^\infty\frac{H_\cap(a)}{a^2}\dd a
       =2\int_0^\infty\frac{e^{-a}-e^{-2a}}a\dd a=2\log2.
\]
Combining this with the $H_2$ and $P_r$ integrals already computed,
we have
\[
 \int_0^\infty\frac{H_2(a)}{a^2}\dd a=1,\quad
 \int_0^\infty\frac{H_\cap(a)}{a^2}\dd a=2\log2,\quad
 \int_0^\infty\frac{P_r(a)}{a^2}\dd a=\frac1{r(r-1)},
\]
which proves all the leading center asymptotics.

\noindent\emph{Positive definiteness.}
Fix $a>0$ and real coefficients $b_0,\ldots,b_K$.
If $\sum_{i=0}^K b_if_i$ has zero variance,
it is constant at every state $(U,V)$, since each state has
positive Poisson probability. At $(0,0)$ its value is zero.
The following states then give, in order,
\begin{align*}
 (U,V)=(1,1)&\quad\Longrightarrow\quad b_1=0,\\
 (U,V)=(K+1,0)&\quad\Longrightarrow\quad b_0=0,\\
 (U,V)=(r,0)&\quad\Longrightarrow\quad b_r=0
                  \quad(2\le r\le K).
\end{align*}
Thus every nonzero coefficient vector has positive local
variance for each $a>0$. Its integral in \eqref{eq:integral}
is positive, so $\Sigma_K$ is strictly positive definite.

\noindent\emph{Covariances of the model counts.}
We now combine the two parts of the covariance calculation.
By the linear transformation defining $Z$ and \eqref{eq:mixedcov},
the sum of the absolute values of covariances between distinct
output primes over $z\le q,\ell\le Q$ is
$o_K(\frac{\log_2x}{\log_3x})$.
Together with \eqref{eq:integral}, this first gives variance
$O_K(\frac{\log_2x}{\log_3x})$ for each sum restricted to $z\le q\le Q$, as well
as the asserted covariance asymptotics for those sums.
The centered omitted sums have $L^2$ norm $o_K(1)$ by
\eqref{eq:smalloutputs}. Cauchy--Schwarz shows that adding back
the small output primes changes any covariance by $o_K(\frac{\log_2x}{\log_3x})$. Hence
\begin{equation}
 \frac{\log_3x}{\log_2x}\Cov_B(Z_i,Z_j)
       =(\Sigma_K)_{ij}+o_K(1).
 \label{eq:modelcov}
\end{equation}

\subsection{The Gaussian limit in the model}
The counts for different output primes are not independent.
We keep the joint distribution at each $q$ and compare the
different output-prime blocks with independent blocks.
The exact-rank coordinates are differences of increasing
thresholds; their signed linear combinations need not be
associated. We therefore make the comparison using the
associated threshold coordinates $X_{h,q}$.

Fix a real coefficient vector $b=(b_0,\ldots,b_K)$. For each
$q$, the definitions of $Z_{i,q}$ give
\[
 \sum_{i=0}^K b_iZ_{i,q}
       =\sum_{h\in\mathcal H}c_hX_{h,q},
\]
where
\begin{align*}
 c_2&=b_0+b_2,\qquad c_\cap=b_1,\qquad c_{K+1}=-b_K,\\
 c_k&=b_k-b_{k-1}\qquad(3\le k\le K).
\end{align*}
The last range is empty when $K=2$. These coefficients are
fixed and may have either sign.

\noindent\emph{Comparison of characteristic functions.}
The standard covariance inequality for associated vectors
\cite[(2) and the observation following Definition 3]{BS}
applies to any two disjoint blocks
$\mathbf U=(U_i)$ and $\mathbf V=(V_j)$ of the threshold
vector. For real functions $F(\mathbf U)$ and $G(\mathbf V)$,
let $a_i$ and $d_j$ be their coordinate Lipschitz bounds:
changing only $U_i$ by $s$ changes $F$ by at most $a_i|s|$,
and the corresponding bound for $G$ is $d_j|s|$. Then
\[
 |\Cov_B(F(\mathbf U),G(\mathbf V))|
       \le\sum_{i,j}a_id_j\Cov_B(U_i,V_j).
\]
For complex-valued functions, apply this bound to the four
pairs of real and imaginary parts. This gives the same bound
with a factor $4$. Here covariance is bilinear, without complex
conjugation.

For $t\in\mathbb R$, we apply this bound successively to the
exponential for the preceding output primes and the exponential
for the next prime. Their coordinate Lipschitz constants are
$|tc_h|$. Summing these differences compares the full
characteristic function with the product below; the other factors
all have modulus at most one. Centering multiplies both terms
by the same phase, so the bound also holds for centered coordinates:

\begin{align}
 &\left|\E_B\exp\left(it\sum_{z\le q\le Q}\sum_{h\in\mathcal H}
                      c_hX_{h,q}\right)\right.\notag\\*
 &\left.\qquad-\prod_{z\le q\le Q}\E_B\exp\left(it\sum_{h\in\mathcal H}
                      c_hX_{h,q}\right)\right|\notag\\
 &\quad\le4t^2\sum_{z\le q<\ell\le Q}
          \sum_{h,j\in\mathcal H}|c_hc_j|
                      \Cov_B(X_{h,q},X_{j,\ell}).
 \label{eq:blockCF}
\end{align}

Apply \eqref{eq:blockCF} to the centered variables, with $t$
replaced by $t/\sqrt{\frac{\log_2x}{\log_3x}}$. For the fixed $b$ and $t$, the
comparison error is bounded by
\begin{align*}
 &\frac{4t^2\log_3x}{\log_2x}\sum_{z\le q<\ell\le Q}
       \sum_{h,j\in\mathcal H}|c_hc_j|
                       \Cov_B(X_{h,q},X_{j,\ell})\\
 &\qquad\ll_{K,b,t}\frac{\log_3x}{\log_2x}
                +\frac{(\log_4x)^2}{\log_3x}
       =o_{K,b,t}(1),
\end{align*}
by \eqref{eq:mixedcov}.
The absolute values of the coefficients allow either sign
in the original linear combination.
The product is the characteristic function of independent
vector blocks with the original joint distribution at each $q$.

\noindent\emph{Limit for the independent blocks.}
It remains to prove the limit for these blocks.
For the fixed vector $b\ne0$, the centered
block sums before division by $\sqrt{\frac{\log_2x}{\log_3x}}$ are uniformly
$O_{K,b}(1)$. The blocks with $q<z$ have total within-block
variance at most $O_{K,b}(z)=o_{K,b}(\log_2x/\log_3x)$.
Thus \eqref{eq:integral} gives the following sum of variances
over $z\le q\le Q$:
\[
 \frac{\log_2x}{\log_3x}\,b^{\mathsf T}\Sigma_K b
       +o_{K,b}\left(\frac{\log_2x}{\log_3x}\right),
\]
where $b^{\mathsf T}$ denotes the transpose of $b$.
Independence makes this the variance of their sum, and its
leading coefficient is positive by positive definiteness.

The standard deviation of the sum therefore tends to infinity,
so every Lindeberg truncation term is zero for sufficiently
large $x$. Lindeberg's theorem gives the Gaussian limit for
the independent blocks. By the characteristic-function comparison
and L\'evy's continuity theorem, the same scalar limit holds in
the original model. We then apply Cram\'er--Wold and restore
the centered small-output sums to obtain
\begin{equation}
 \sqrt{\frac{\log_3x}{\log_2x}}
       (Z_i-\E_BZ_i)_{i=0}^K
       \Longrightarrow N_{K+1}(0,\Sigma_K).
 \label{eq:modelCLT}
\end{equation}

\subsection{Return to the arithmetic counts}
We have proved the joint limit and covariance asymptotics in
the Bernoulli model. To obtain the same conclusions for
integers $n\le x$, we first compare the truncated arithmetic
counts with this model. We then add back the omitted primes.
Write $T=(T_0,\ldots,T_K)$, where
\[
 T_0=W,\qquad T_1=C,\qquad T_r=E_r\quad(2\le r\le K).
\]
Define their truncated approximations by
\begin{align*}
 \widehat T_0&=\mathcal C_2(n;y,Q),\\
 \widehat T_1&=\mathcal C_\cap(n;y,Q),\\
 \widehat T_r&=\mathcal C_r(n;y,Q)-\mathcal C_{r+1}(n;y,Q)
       \quad(2\le r\le K).
\end{align*}
Write $\widehat T=(\widehat T_0,\ldots,\widehat T_K)$.
Every coordinate of $\widehat T$ and $Z$ is between zero and
$Q$, including the exact-rank differences. We apply the Kubilius
comparison \eqref{eq:TVinput} to the entire vector. The
deterministic maps described above then give
\begin{align*}
 d_{\rm TV}(\widehat T,Z)&\le\delta,\\
 |\E_x\widehat T_i-\E_BZ_i|&\ll Q\delta,\\
 |\E_x\widehat T_i\widehat T_j-\E_BZ_iZ_j|
       &\ll Q^2\delta.
\end{align*}
Since all these coordinates are bounded by $Q$, the difference
between the products of their means is also $O(Q^2\delta)$.
Thus the covariance comparison has the same error bound.
Equations \eqref{eq:TVsmall} and \eqref{eq:Zmean} give the mean
estimates stated in the theorem for the corresponding truncated
coordinates. Together with \eqref{eq:modelcov}, they give
\begin{equation}
 \frac{\log_3x}{\log_2x}
       \Cov_x(\widehat T_i,\widehat T_j)
       =(\Sigma_K)_{ij}+o_K(1).
 \label{eq:truncatedmoments}
\end{equation}
Applying the same centering and scaling to both vectors
preserves the total variation bound. Thus the normalized
vectors centered at $\E_BZ$ have asymptotically the same law.
The difference between $\E_x\widehat T$ and $\E_BZ$, after
multiplication by $\sqrt{\log_3x/\log_2x}$, tends to zero in
each coordinate, since $Q\delta=o(\sqrt{\log_2x/\log_3x})$.
Consequently, \eqref{eq:modelCLT} also holds for the truncated
arithmetic vector centered at its own mean.

We next remove the truncation. For the exact-rank coordinates,
write the error as a group comparison and two truncation errors:
\begin{align*}
 E_r-\widehat T_r
 &=[E_r-(\mathcal C_r-\mathcal C_{r+1})]\\
 &\quad+[\mathcal C_r-\mathcal C_r(n;y,Q)]\\
 &\quad-[\mathcal C_{r+1}-\mathcal C_{r+1}(n;y,Q)].
\end{align*}
The coordinates $W,C$ have the same decomposition with one
occupancy error. Lemma \ref{lem:groups} controls the group
comparison, and Lemma \ref{lem:trunc} controls each occupancy
error. For the first moments, use the triangle inequality and
Cauchy--Schwarz on the group comparison. For the second
moments, use $(a+b+c)^2\le3(a^2+b^2+c^2)$. This gives,
for every coordinate,
\begin{align}
 \E_x|T_i-\widehat T_i|
       &\ll_K\log_3x,\label{eq:finalLone}\\
 \E_x|T_i-\widehat T_i|^2
       &\ll_K(\log_3x)^2+(1+\log_4x)^4\notag\\
       &\quad+\left(\frac{\log_2x}{\log_3x}\right)^{2/3}
          =o_K\left(\frac{\log_2x}{\log_3x}\right).
       \label{eq:finalLtwo}
\end{align}
The variance of each error is at most its second moment.
Chebyshev's inequality therefore shows that each normalized
centered error tends to zero in probability. Since $K$ is fixed,
the same is true for the whole error vector. Slutsky's theorem
then gives the joint limit for
$(T_i-\E_xT_i)_{i=0}^K/\sqrt{\frac{\log_2x}{\log_3x}}$.
Combining \eqref{eq:finalLone} with the truncated mean estimates
proves the arithmetic mean estimates in the theorem.
Since $\log_3x=o(\sqrt{\frac{\log_2x}{\log_3x}})$, those estimates allow replacement
of the exact means by the explicit centers in \eqref{eq:centers},
giving \eqref{eq:joint}.

It remains to prove the covariance formula for the full
arithmetic counts. Put $\varepsilon_i=T_i-\widehat T_i$.
The truncated standard deviations
are $O_K(\sqrt{\frac{\log_2x}{\log_3x}})$ by \eqref{eq:truncatedmoments}, while
$\|\varepsilon_i\|_2=o_K(\sqrt{\frac{\log_2x}{\log_3x}})$. Expanding the difference of
covariances gives
\begin{align*}
 &\Cov_x(T_i,T_j)-\Cov_x(\widehat T_i,\widehat T_j)\\*
 &\quad=\Cov_x(\varepsilon_i,\widehat T_j)
       +\Cov_x(\widehat T_i,\varepsilon_j)
       +\Cov_x(\varepsilon_i,\varepsilon_j)\\
 &\quad=o_K\left(\frac{\log_2x}{\log_3x}\right).
\end{align*}
Cauchy--Schwarz and $\Var_x\varepsilon_i\le\E_x\varepsilon_i^2$
bound the absolute value of this difference by
\begin{align*}
 &\|\varepsilon_i\|_2
       \sqrt{\Var_x\widehat T_j}
 +\|\varepsilon_j\|_2
       \sqrt{\Var_x\widehat T_i}
 \\
 &\quad+\|\varepsilon_i\|_2\|\varepsilon_j\|_2
       =o_K\left(\frac{\log_2x}{\log_3x}\right).
\end{align*}
This proves the arithmetic covariance assertion of the theorem.
The bounded mean-square $\sigma$ comparison in
Lemma \ref{lem:groups} also proves Corollary \ref{cor:sigma}.
This completes the proof of Theorem \ref{thm:joint}. \qed

\begin{remark}
The rank bound $K$ is fixed throughout the proof. We do not obtain a
growing-rank or functional central limit theorem, nor a rate
of convergence in \eqref{eq:joint}. The mean-square estimates
are needed to return to the arithmetic covariances: an $L^1$
truncation would suffice for the Gaussian law but would not
by itself justify the arithmetic covariance asymptotics.
\end{remark}

\section*{Declaration of AI use}
The author used ChatGPT and Codex (OpenAI) for literature searches,
manuscript preparation, and the development and checking of mathematical
arguments and calculations. The author takes full responsibility for
the manuscript.

\end{document}